\documentclass{article}
\usepackage{amsthm, amssymb, latexsym}
\usepackage[cp1251]{inputenc}
\usepackage[english,russian,ukrainian]{babel}
\usepackage[T2A]{fontenc}

\begin{document}
\thispagestyle{empty}
\setcounter{page}{1}

\newtheorem{theorem}{Theorem}
\newtheorem{lemma}{Lemma}
\newtheorem{corollary}{Corollary}
\newtheorem{proposition}{Proposition}
\newtheorem{definition}{Definition}
\newtheorem{assertion}{Assertion}
\newcommand{\rank}{\mbox{rank}}
\newcommand{\Aut}{\mbox{Aut }}
\newcommand{\Ker}{\mbox{Ker }}
\renewcommand{\refname}{}
\theoremstyle{plain}
\mathsurround 2pt

\begin{center} \renewcommand{\baselinestretch}{1.3}\bf {ON FACTOR GROUPS OF SOME GROUPS} \end{center}
\begin{center} \textsc {\emph{N.S. CHERNIKOV}}
\end{center}
\begin{center} {Institute of Mathematics of NASU,\\
3 Tereshenkivska Str., Kyiv 01601, Ukraine}
\end{center}

\renewcommand{\abstractname}{}

\begin{abstract}
Let for a prime $p$, $\mathfrak{X}$ (respectively $\mathfrak{Y}$) be the class of all $p$-biprimitively finite (respectively periodic $p$-conjugatively biprimitively finite) groups and $G\in \mathfrak{X}$ (respectively $G\in \mathfrak{Y}$), $V$ be a periodic subgroup of $G$ having an ascending series of normal in $G$ subgroups such that each its factor is an almost layer-finite group or a locally graded group of finite special rank, or a $WF$-group with $min-q$ on all primes $q$. We prove that $G/V \in \mathfrak{X}$ (respectively $G/V\in \mathfrak{Y}$). Also some interesting and useful preliminary results are obtained.
\end{abstract}

Below $p$ and $q$ are primes, $\mathbb{N}$ is the set of all natural numbers, $\mathbb{P}$ is the set of all primes; $min$ (respectively $min-q$) is the minimal condition on ($q$-) subgroups, $\times$ is the sign of the direct product; $H~sn~G$ denotes that $H$ is a subnormal subgroup of the group $G$. Other notations in the present paper are standard. The empty set is not reckoned finite. Recall that for a class $\mathfrak{X}$ of groups, an almost $\mathfrak{X}$-group is a finite extension of a group belonging to $\mathfrak{X}$.\

Remind definitions of some classes of groups: the group $G$ is called $p$-(conjugatively) biprimitively finite, if for every its finite subgroup $K$, each subgroup of $N_G(K)/K$, generated by two (conjugate) elements of order $p$, is finite; a group is called (conjugatively) biprimitively finite, if for any $p$, it is $p$-(conjugatively) biprimitively finite (V.P.Shunkov, 1970 - 1973). Conjugatively biprimitively finite groups are also called Shunkov groups (V.D.Mazurov, 2003). These classes of groups are very large and contain, for instance, all locally finite groups, all $2$-groups. Many deep results of Shunkov's School are connected with groups of these classes (see, for instance, \cite{SSh}).\

Futher, a group is called layer-finite, if the set of all its elements of each order is finite or empty (S.N.Chernikov, see, for instance, \cite{CHS80}). In consequence of Dietzmann's Lemma (see, for instance, \cite{Rob}), a layer-finite group is locally (finite and normal). A layer-finite group is also named an $FL$-group \cite{Rob}.\

Remind that the $FC$-centre of the group $G$ is the following its subgroup: $\{g\in G: |G:C_G(g)|<\infty\}$, and the $FC$-hypercentre of $G$ is the limit of its ascending normal series defined by the rules: $G_0=1$, $G_{\alpha+1}/G_\alpha=$ the $FC$-centre of $G/G_\alpha$, and $G_\lambda=\cup_{\beta<\lambda}G_\beta$, where $\alpha$ is an ordinal and $\lambda\neq 0$ is a limit ordinal (see, for instance, \cite{Rob}). The $FC$-centre (resp. $FC$-hypercentre) of $G$ contains its centre (resp. hypercentre).\

Remind: a group possessing an ascending series with finite and locally nilpotent factors is called a $WF$-group (B.I.Plotkin, \cite{Pl}). The class of all $WF$-groups is large. It contains, for instance, all $RN^*$-groups, all radical in the sense of B.I.Plotkin groups, all groups possessing an ascending series with $FC$-factors (at the same time, all hyperfinite groups).\

Remind: a group such that any its finitely generated subgroup $\neq 1$ possesses a subgroup of finite index $\neq 1$ is called locally graded (S.N.Chernikov, see, for instance, \cite {CHS80}). The class of all locally graded groups is very large. It contains, for example, all locally finite, locally solvable, residually finite, linear, radical in the sense of B.I.Plotkin groups, all $WF$-groups, all $RN$-groups and, at the same time, all groups of all Kurosh-Chernikov classes \cite {KarMer}.\

The main result of the present paper is the following theorem.\

\textbf{Theorem.}\emph{
Let $G$ be a $p$-biprimitively finite (respectively periodic $p$-conjugatively biprimitively finite) group and $V$ its normal periodic subgroup. Suppose that for $V$ even if one of the following conditions is fulfilled:
\begin{itemize}
 \item [(a)] It possesses an ascending series ${\mathcal{M}}$ of normal in $G$ subgroups every factor of which is an almost layer-finite group or a locally graded group of finite special rank or a $WF$-group with $min-q$ on all $q$ (i.e. with Artinian $q$-subgroups on all $q$).
 \item [(b)] It possesses an ascending series with finite factors of normal in $G$ subgroups.
 \item [(c)] It belongs to the $FC$-hypercentre of $G$.
\end{itemize}
Then $G/V$ is $p$-biprimitively finite (respectively $p$-conjugatively biprimitively finite).}

The following corollaries are immediate consequences of Theorem.
\begin{corollary}\label{2} Let $G$ be a $p$-biprimitively finite (respectively periodic $p$-conjugatively biprimitively finite) group and $V$ its normal periodic subgroup. If $V$ is almost layer-finite or locally graded of finite special rank, or a $WF$-group with $min-q$ on all $q$, then $G/V$ is $p$-biprimitively finite (respectively $p$-conjugatively biprimitively finite).
\end{corollary}
\begin{corollary}\label{3} Let $G$ be a biprimitively finite (respectively periodic Shunkov) group and $V$ its normal periodic subgroup. If $V$ is almost layer-finite or locally graded of finite special rank, or a $WF$-group with $min-q$ on all $q$, then $G/V$ is biprimitively finite (respectively Shunkov).
\end{corollary}
\begin{corollary}\label{4} Let $G$ be a biprimitively finite (respectively periodic Shunkov) group and $V$ its subgroup possessing an ascending series with finite factors of normal in $G$ subgroups. Then $G/V$ is biprimitively finite (respectively Shunkov).
\end{corollary}
Preface the proof of Theorem with the following propositions.
\begin{proposition}\label{5} Let $G$ be a group having an ascending normal series $G_0=1\subset ... \subset G_\gamma=G$ with finite factors, $M$ a finite set of elements of $G$, $\Omega_1, ..., \Omega_n$ operator groups of $G$, generated by finite sets of elements of finite order, such that for any $k=1, ..., n$ and ordinal $\alpha<\gamma$, ${G_\alpha}^{\Omega_k}=G_\alpha$. Then:
\begin{itemize}
 \item [(i)] All finite subgroups $F$ of $G$, for which $F^{\Omega_k}=F$, $k=1, ..., n$, constitute a local system of $G$.
 \item [(ii)] For any $\Omega_k$, $|\Omega_k:C_{\Omega_k}(M))|<\infty$, $\Omega_k/C_{\Omega_k}(G)$ is residually finite, and in the case when $\Omega_k$ is a group of automorphisms of $G$, $\Omega_k$ is residually finite.
\end{itemize}
\end{proposition}
\begin{proposition}\label{6} Let $G$ be a group having a layer-finite subgroup $L$ of finite index. Then G is a locally finite group with Chernikov $q$-subgroups on all $q$, its $FC$-centre $F$ is layer-finite  and $|G:F|<\infty$. Further, an arbitrary $FC$-subgroup $B$ of $G$ (in particular, an arbitrary layer-finite subgroup B of $G$) belongs to $F$, if $|G:B|<\infty$ or $B~sn~G$.
\end{proposition}
\begin{proposition}\label{7} In each of the following cases the periodic group $G$ possesses an ascending series with finite factors of characteristic subgroups:
\begin{itemize}
 \item [(a)] $G$ is almost layer-finite.
 \item [(b)] $G$ is locally graded of finite special rank.
 \item [(c)] $G$ is a $WF$-group with $min-q$ on all $q$ (i.e. with Artinian $q$-subgroups on all $q$).
\end{itemize}
\end{proposition}
\begin{proposition}\label{8} Let $G$ be a group, $H$ its subgroup generated by a finite set of elements of finite order and $V$ its normal periodic subgroup. Suppose that for $V$ even if one of the conditions (a)-(c) of Theorem is fulfilled. Then:
\begin{itemize}
 \item [(i)] All finite subgroup of $V$, normalized by $H$, constitute a local system of $V$.
 \item [(ii)] $H/C_H(V)$ is residually finite.
\end{itemize}
\end{proposition}
The following proposition is an immediate consequence of Proposition \ref{8}.
\begin{corollary}\label{9} Let $G$ be a group, $H$ its subgroup generated by a finite set of elements of finite order and $V$ its normal periodic subgroup.Suppose that $V$ is almost layer-finite or locally graded of finite special rank, or a $WF$-group with $min-q$ on all $q$. Then the statements of Proposition \ref{8} are valid.
\end{corollary}

\textbf{\emph{Proof of Proposition 1.}}\

Let (i)  be false and $G$ be a counter-example with minimal $\gamma$. Then for some finite set $X$ of elements of $G$, there are no finite subgroups $F \supseteq X$ such that $F^{\Omega_k}=F$, $k=1, ..., n$.\

For some limit infinite ordinal $\nu\leq\gamma$, $|G:G_\nu|<\infty$. Let $T$ be a transversal to $G_\nu$ in $G$, $\mathfrak{M}_k$ be a finite set of finite cyclic subgroups of $\Omega_k$ that generate $\Omega_k$, and $\bigcup\limits_{k=1}^n \mathfrak{M}_k=\{\Delta_1, \ldots, \Delta_m\}$; $H_j=\langle X^{\Delta_j}, T^{\Delta_j}\rangle$, $j=1,\ldots,m$, and $H=\langle H_1,\ldots,H_m\rangle$. Then $H_j^{\Delta_j}=H_j$ and $H=H_j(H\cap G_\nu)$, $j=1,\ldots,m$.

In view of Corollary [3,P.35] (for instance), $G$ is locally finite. Therefore because of $H$ is finitely generated, it is finite. Since $\nu$ is infinite limit and $H\cap G_\nu$ is finite, for some ordinal $\beta<\nu$, $H\cap G_\nu\subseteq G_\beta$. Put $\zeta=\beta$, if $H\subseteq G_\beta$, and $\zeta=\beta+1$, if $H\nsubseteq G_\beta$; $L=HG_\beta$ and $L_\alpha=G_\alpha$, $\alpha<\zeta$. Then $X\subseteq L$, $\zeta<\gamma$, and $L_0=1\subset ... \subset L_\zeta=L$ is an ascending normal series of $L$ with finite factors and also for $\alpha<\zeta$, ${L_\alpha}^{\Omega_k}=L_\alpha$, $k=1, ..., n$. Further,
$$L^{\Delta_j}=(HG_\beta)^{\Delta_j}=(H_j(H\cap G_\nu)G_\beta)^{\Delta_j}=(H_jG_\beta)^{\Delta_j}=$$
$$={H_j}^{\Delta_j}{G_\beta^{\Delta_j}}=H_jG_\beta=H_j(H\cap G_\nu)G_\beta=HG_\beta=L, j=1,\ldots,m.$$
Therefore because of $\Omega_k$ is generated by some subgroups $\Delta_j$, $L^{\Omega_k}=L$. Thus ${L_\alpha}^{\Omega_k}=L_\alpha$, $\alpha\leq \zeta$ and $k=1,\ldots,n$. Further, since $\zeta<\gamma$, $L$ is not a counter - example to (i). Therefore for some finite subgroup $F$ of $L$, $X\subseteq F$ and $F^{\Omega_k}=F$, $k=1,\ldots,n$, which is a contradiction.
Thus (i) is correct.\

Further, in view of (i), for some finite $F\subseteq G$, $M\subseteq F=F^{\Omega_k}$. So $|\Omega_k:C_{\Omega_k}(M)|\leq |\Omega_k:C_{\Omega_k}(F)|<\infty$. Therefore because of $C_{\Omega_k}(G)$ is the intersection of subgroups $C_{\Omega_k}(M)$ by all finite $M\subseteq G$, $\Omega_k/C_{\Omega_k}(G)$ is residually finite. In the last case, $C_{\Omega_k}(G)=1$ and so $\Omega_k$ is residually finite.\

Proposition is proven.\

\textbf{\emph{Proof of Proposition 2.}}\

First, if $|G:B|<\infty$, then for $b\in B$, $|G:C_G(B)|<\infty$. So $B\subseteq F$.\

Since $|G:L|<\infty$ and $L$ is an $FC$-group, $L\subseteq F$. Consequently, $|G:F|<\infty$.\

Let $X$ be any finitely generated subgroup of $G$. Then $|X:X\cap L|<\infty$. In consequence of Schreier's Theorem (see, for instance, [6,P.228]), $X\cap L$ is finitely generated. Since $L$ is locally finite, $X\cap L$ is finite. So $X$ is finite. Thus $G$ is locally finite.\

Let $Q$ be a $q$-subgroup of $G$. By Theorem 3.2 \cite{CHS80}, $Q\cap L$ is Chernikov. Since $|Q:Q\cap L|<\infty$, obviously $Q$ is Chernikov too.\

Let the statement of the present proposition, relating to the case when $B~sn~G$, be incorrect, and $G$ be a counter-example to this statement such that the length of some passing through $B$ finite series $G_0=1\subset B=G_1\subset \ldots \subset G_n=G$ of $G$ is minimal. For the layer-finite subgroup $L\cap G_{n-1}$, $|G_{n-1}:L\cap G_{n-1}|<\infty$. Then the $FC$-centre $K$ of $G_{n-1}$ contains $B$. Since obviously $K$ is locally (finite and normal) and satisfies for each $q$ $min-q$, it is layer-finite (Theorem 3.7 \cite{CHS80}). Let $a\in K$. Then for any $g\in G$, $|<a>|=|<a^g>|$ and $a^g\in K$. Therefore because $K$ is layer-finite, the set $\{a^g: g\in G\}$ is finite. So $a\in F$. Thus $B\subseteq K\subseteq F$, which is a contradiction.\

Proposition is proven.\

\textbf{\emph{Proof of Proposition 3.}}\

(a) In view of Proposition 2, $G$ is locally finite and contains some characteristic layer-finite subgroup $F$ of finite index. For any $k\in \mathbb{N}$, $\langle g\in F: |\langle g\rangle|\leq K\rangle$ is obviously a finite characteristic subgroup of G. Clearly, all distinct subgroups among subgroups: $\langle g\in F: |\langle g\rangle|\leq K\rangle$, $k=1,2,\ldots$, $F$, $G$, constitute a required series of $G$.

Let $N$ be the product of all subgroups that possess an ascending series with finite factors of characteristic subgroups of $G$. Then $N$ has the same series, and $G/N$ has no finite characteristic subgroups $\neq 1$. Consequently: $G$ has the same series iff $G=N$; if $G/N\neq 1$, then $G/N$ is infinite.\

Let $B$ be a characteristic Chernikov subgroup of $G/N$. If $B\neq 1$, then $B$ contains some subgroup $H$ of finite index, which is a direct product of finitely many quasicyclic subgroups. Let $p$ be a prime for which $H$ has an element of order $p$. Then $\{g\in H: |\langle g\rangle|\leq p\}\neq 1$. But $\{g\in H: |\langle g\rangle|\leq p\}$ is clearly a finite characteristic subgroup of $G/N$, which is a contradiction. Thus $G/N$ has no characteristic Chernikov subgroups $\neq 1$.\

(b) Suppose $G/N\neq 1$. In view of Theorem [7], $G$ and, at the same time, $G/N$ are locally finite and almost hyperabelian. Let $K$ be a normal hyperabelian subgroup of finite index of $G/N$. Then $K\neq 1$. So $K$ contains some normal abelian subgroup $A\neq 1$. For some $p$, $O_p(A)\neq 1$. Since $O_p(A)\unlhd K$ and $O_p(K)\unlhd G/N$, $O_p(G/N)\neq 1$. Since $O_p(G/N)$ is locally finite of finite special rank, it is Chernikov \cite{Myag}, which is a contradiction. Thus $G=N$, so (b) is correct.\

(c) Suppose $G/N\neq 1$. In view of Proposition 1.1 \cite{CHS80} and Corollary [3,P.35] (for instance), $G$ and, at the same time, $G/N$ are locally finite. Therefore for each $q$ all $q$-subgroups of $G$ are Chernikov (Theorem 1.5 \cite {CHS80}). Then by virtue of Theorem 3.13 \cite{KegWeh}, for each $q$ all $q$-subgroups of $G/N$ are Chernikov. Since $O_q(G/N)$ is Chernikov, $O_q(G/N)=1$ (see above).

Since, obviously, $G/N$ is a $WF$-group, it has some ascendent subgroup $L\neq 1$ which is locally nilpotent or finite.\

In the first case, $L=\times_{q\in \mathbb{P}} L_q$ where $L_q$ are Sylow $q$-subgroups of $L$. For some $p$, $L_p\neq 1$. By Lemma 2.1 \cite{CHS59}, $O_p(G/N)\neq 1$, which is a contradiction.\

Thus $L$ is finite.

Let $M$ be a subnormal subgroup of $L$ of minimal $\neq 1$ order. Then $M$ is an ascendant finite simple subgroup of $G/N$. It is non-abelian (see above). By Lemma 2.1 \cite{CHS59}, $G/N$ contains some characteristic subgroup $R\neq 1$ which is a direct product of subgroups isomorphic to $M$. Then $R$ is infinite. So because $M$ is finite, $R$ contains an infinite subgroup which is a direct product of subgroups of order $p$ for some $p$. But this subgroup is not Chernikov, which is a contradiction.\

Thus $G=N$. So (c) is correct.\

Proposition is proven.\

\textbf{\emph{Proof of Proposition 4.}}\

In view of Proposition 3, $\mathcal{M}$ is contained in an ascending series of $V$ with finite factors of normal in $G$ subgroups, i.e. (a) implies (b).\

Suppose (c) is fulfilled. Then it is easy to see: $V$ has an ascending series of normal in $G$ subgroups such that every its factor $A/B$ contains some element $g$, for which $|G/B:C_{G/B}(g)|<\infty$ and $A/B=\langle g^{G/B}\rangle$. In view of Dietzmann's Lemma, $\langle g^{G/B}\rangle$ is finite. So (c) implies (b).\

Thus (b) is necessarily fulfilled.\

Setting in Proposition 1 $G=V$, $n=1$ and $\Omega_1=H$ and applying this proposition, conclude that the present proposition is true.\

\textbf{\emph{Proof of Theorem.}}\

Let $L/V$ be a finite subgroup of $G/V$, $U/V$ the normalizer of $L/V$ in $G/V$ and $(T/V)/(L/V)$ a subgroup of $(U/V)/(L/V)$, generated by two (conjugate) elements of order $p$. It is easy to see that for some $a, b\in N_G(L)\setminus L$, $a^p, b^p\in L$ and $T=\langle a, b\rangle L$ (resp. for some $a\in N_G(L)\setminus L$ and $g\in N_G(L)$, $a^p\in L$ and $T=\langle a, a^g\rangle L$). Note that $a, b, g$ are elements of finite order.\

It is easy to see that for $L$ even if one of the conditions (a)-(c) is fulfilled. Therefore by virtue of Proposition 4, for some finite subgroup $F$ of $L$, $\langle a, b\rangle\subseteq N_G(F)$ and $a^p, b^p\in F$ (respectively $\langle a, g\rangle\subseteq N_G(F)$ and $a^p\in F$ and, at the same time, $a^{pg}\in F$). Make more precise: in Proposition 4, we set $H=\langle a, b\rangle$ or $H=\langle a, g\rangle$ respectively. Since $aF$ and $bF$ (resp. $aF$ and $a^gF$) are (conjugate) elements of order $p$ of $N_G(F)/F$, the subgroup $\langle aF, bF\rangle$ (respectively $\langle aF, a^gF\rangle$) of $N_G(F)/F$ is finite. Therefore $\langle a, b\rangle$ (respectively $\langle a, a^g\rangle$) is finite. At the same time, $(T/V)/(L/V)$ is finite. Thus, the present theorem is true.\

The following assertion is an immediate consequence of Proposition \ref{5}.\

\textbf{Assertion.} \emph{Let $G$ be a Chernikov group, $M$ a finite set of elements of $G$, $\Omega_1, ..., \Omega_n$ operator groups of $G$, generated by finite sets of elements of finite order. Then the statements (i) and (ii) of Proposition 1 are valid.}


\begin{thebibliography}{99}
\bibitem{SSh} \textit{Senashov V.I., Shunkov V.P.} Groups with finiteness conditions. - Novosibirsk: Publ. House Siberian Branch Russian Acad. Sci., 2001. - 336 p. (in Russian)
\bibitem{CHS80} \textit{Chernikov S.N.} Groups with prescribed properties of the system of subgroups. -- Мoskva: Nauka, 1980. -- 384~p.(in Russian)
\bibitem{Rob} \textit{Robinson D.J.S.}Finiteness conditions and generalized soluble groups. Pt 1. -- Berlin -- Heidelberg -- New York: Springer, 1972. -- 210~p.
\bibitem{Pl} \textit{Plotkin B.I.}Radical and semi-simple group// Trudy Moskov. Mat. Obsh. -- 1957. -- \textbf{6}. -- P.~299--336 (in Russian).
\bibitem{KarMer} \textit{Kargapolov M.I., Merzljakov Yu.I.}Fundamentals of the theory of groups. -- Мoskva: Nauka, 1972. -- 240~p.(in Russian)
\bibitem{Kur} \textit{Kurosh A.G.}The theory of groups. 3-d ed. -- Мoskva: Nauka, 1967. -- 648~p. (in Russian)
\bibitem{CHN} \textit{Chernikov N.S.}One theorem on groups of finite special rank// Ukr.Mat.Z. -- 1990. -- \textbf{42}, №7. - P.~962--970 (in Russian).
\bibitem{Myag} \textit{Myagkova N.N.}On groups of finite rank// Izv. Acad.Nauk SSSR. Ser. Mat. -- 1949. -- \textbf{13}, №6. - P.~495--512 (in Russian).
\bibitem{KegWeh} \textit{Kegel O.H., Wehrfritz B.A.F.}Locally finite groups. -- Amsterdam -- London: North-Holland Publ. Co, 1973. -- 210~p.
\bibitem{CHS59} \textit{Chernikov S.N.}Finiteness conditions in the general theory of groups// Uspehi. Mat. Nauk. -- 1959. -- \textbf{14}, №5. -- P.~45-96 (in Russian).
\end{thebibliography}
\end{document}